\documentclass{article}
\begin{document}
\centerline{\textbf{\large A General Method for Constructing  Ramanujan Formulas for $1/\pi^{\nu}$}}
\[
\]
\centerline{\bf Nikos Bagis}\\
\centerline{e-mail: nikosbagis@hotmail.gr}
\[
\]
\begin{quote}
\begin{abstract}
In this article we give the theoretical background for generating Ramanujan type $1/\pi^{2\nu}$ formulas. As applications of our method we give a general  construction of $1/\pi^4$ series and examples of $1/\pi^6$ series. We also study the elliptic alpha function whose values are useful for such evaluations.   
\end{abstract}

\bf keywords \rm{$\pi$-formulas; elliptic alpha function; Ramanujan; elliptic functions; singular modulus; approximations; numerical methods}
\end{quote}

\section{Introduction}

 The standard definitions of the Elliptic Integrals of the first and second kind respectively (see [9],[4],[18]) are:
\begin{equation}
K(x)=\int^{\pi/2}_{0}\frac{dt}{\sqrt{1-x^2\sin^2(t)}} \textrm{ and } E(x)=\int^{\pi/2}_{0}\sqrt{1-x^2\sin^2(t)}dt . 
\end{equation}
\begin{equation}
K(x)=\frac{\pi}{2}{}_2F_{1}\left(1/2,1/2,1;x^2\right)=\frac{\pi}{2}\sum^{\infty}_{n=0}\frac{\left(\frac{1}{2}\right)_n\left(\frac{1}{2}\right)_n}{(1)_n}\frac{x^{2n}}{n!}
\end{equation}
In the notation of Mathematica we have
\begin{equation}
K(x)=\textrm{EllipticK}[x^2]\textrm{ and } E(x)=\textrm{EllipticE}[x^2] . 
\end{equation}
Also we have (see [7],[9]):  
\begin{equation}
\dot{K}(k)=\frac{dK(k)}{dk}=\frac{E(k)}{k(1-k^2)}-\frac{K(k)}{k} . 
\end{equation}
The elliptic singular moduli $k_r$ is defined to be the solution of the equation: 
\begin{equation}
\frac{K\left(\sqrt{1-w^2}\right)}{K(w)}=\sqrt{r} .
\end{equation}
In Mathematica notation
\begin{equation}
w=k=k_r=k[r]=\textrm{InverseEllipticNomeQ}[e^{-\pi \sqrt{r}}]^{1/2} . 
\end{equation}
The complementary modulus is given by $k^{'2}_{r}=1-k_r^2$. (For evaluations of $k_r$  see [7],[17],[27])\\  
Also we will need the following relation satisfied by the elliptic alpha function (see [7]):
\begin{equation}
a(r)=\frac{\pi}{4K(k_r)^2}-\sqrt{r}\left(\frac{E(k_r)}{K(k_r)}-1\right) . 
\end{equation}
Our method requires finding derivatives of powers of the elliptic integral $K$  which always can be expressed in terms of $K$ and $E$. Hence from (7) and (4) we need to know $k_r$ and $a(r)$. 
(Here we carry out these evaluations using Mathematica programs.\\ In the literature (see [7] and Wolfram pages) the function $a(r)$ is not  widely known.\\ 
As for  the singular moduli, the elliptic alpha function can evaluated from modular equations. The case of $a(4r)$ is given in [7] chapter 5:
\begin{equation}
a(4r)=(1+k_r)^2a(r)-2\sqrt{r}k_r
\end{equation}   
The authors of [22]  have proven the formula for $a(9r)$:
\begin{equation}
\frac{a(9r)}{\sqrt{r}}-k^2_{9r}=1-\frac{k_{9r}k_{r}}{3M_{3r}}-\frac{k'_{9r}k'_{r}}{3M_{3r}}-\frac{1}{3M_{3r}}-\frac{1}{3M^2_{3r}}+\frac{1}{M^2_{3r}}\left(\frac{a(r)}{\sqrt{r}}-\frac{k^2_{r}}{3}\right)
\end{equation}   
where $M_{3r}$ is a root of the polynomial equation
\begin{equation}
27M^4_{3r}-18M^2_{3r}-8(1-2k^2_r)M_{3r}-1=0
\end{equation}

In Section 2 we review and extend the method for constructing $\pi^{-n}$ series based on $a(r)$. In the next section,
 following Ramanujan (see [4]),  we give and prove a formula for the evaluation of $a(25r)$, when $a(r)$ is known. We conclude this note with some examples.

\section{The general method and the construction of  $1/\pi^4$ and $1/\pi^6$ formulas}

In Mathematica notation (see [28]): 
\begin{equation}
\phi(z)={}_3F_{2}\left(\frac{1}{2},\frac{1}{2},\frac{1}{2};{1,1};z\right)=\frac{4K^2\left(\frac{1}{2}(1-\sqrt{1-z})\right)}{\pi^2} 
\end{equation}
Let
$$
C_{n,r}=\frac{n!}{(n-r)!r!}
$$
and
\begin{equation} 
F_{p}(x):=\phi^{p}(x) 
\end{equation}
and define $c_{p}(n)$, with $p:=2\nu$, $\nu\in\bf N\rm$, by
$$
F_{p}(x)=\left(\sum^{\infty}_{n=0}\frac{C^3_{2n,n}}{64^n} x^n\right)^{p}=\sum^{\infty}_{n=0}c_{p}(n)x^n\eqno :def
$$
Also consider the following equation satisfied by the function $F_{2\nu}(x)$:
$$
A_{p}z^pF^{(p)}_{p}(z)+A_{p-1}z^{p-1}F^{(p-1)}(z)+\ldots+A_{1}zF^{(1)}(z)+A_0F^{(0)}(z)=
$$
\begin{equation}
=\sum^{\infty}_{n=0}c_{2\nu}(n)z^n\left(C_{2\nu}n^{2\nu}+C_{2\nu-1}n^{2\nu-1}+\ldots+C_{1}n^1+C_0\right)=\frac{g}{\pi^{2\nu}}
\end{equation}
where $p=2\nu-1$.
\[
\]
Set $z=1-\left(1-2w\right)^2$. Then if $1-2w=k_r^2$ the quantity $g$ is a function of $k_r$ and $a(r)$ and hence $g$ is algebraic number when $r\in \bf N\rm$ and  $A_{2\nu}$, $A_{2\nu-1}$,\ldots, $A_1$, $A_0$, $C_{2\nu}$, $C_{2\nu-1}$, $\ldots$ $C_1$, $C_0$  can determined from equation (13).\\   
\\
\centerline{\textbf{\large A series representing $1/\pi^4$.}}
\\
For $\nu=2$ we have 
\[
\]
$$
A_4=r^2 (1-2 w)^4 [105a(r)^4-420 a(r)^3 \sqrt{r} w+90 a(r)^2 r w (-1+8 w)-
$$
$$
-20 a(r) r^{3/2} w (1-12 w+32+
+w^2)+r^2 w (-2+43 w-192 w^2+256 w^3)]^{-1}
$$
\[
\]
$$
A_3=[2 r^{3/2} (1-2 w)^2 (a(r) (5-10 w)+\sqrt{r} (3-23 w+28 w^2))]\times
$$
$$
\times[105 a(r)^4-420 a(r)^3 \sqrt{r} w+90 a(r)^2 r w (-1+8 w)-20 a(r) r^{3/2} w (1-12 w+32 w^2)+
$$
$$
+r^2 w (-2+43 w-192 w^2+256 w^3)]^{-1}
$$
\[
\]
$$
A_2=[r (45 a(r)^2 (1-2 w)^2-30 a(r) \sqrt{r} (-1+11 w-30 w^2+24 w^3)+
$$
$$
+r (7-140 w+735 w^2-1400 w^3+880 w^4))][105 a(r)^4-420 a(r)^3 \sqrt{r} w+
$$
$$
+90 a(r)^2 r w (-1+8 w)-20 a(r) r^{3/2} w (1-12 w+32 w^2)+r^2 w (-2+43 w-
$$
$$
-192 w^2+256 w^3)]^{-1}
$$
\[
\]
$$
A_1=[\sqrt{r} (-210 a(r)^3 (-1+2 w)+90 a(r)^2 \sqrt{r} (1-13 w+20 w^2)-
$$
$$
-10 a(r) r (-2+57 w-276 w^2+304 w^3)+r^{3/2} (2-115 w+995 w^2-
$$
$$
-2640 w^3+2080 w^4))][2 (105 a(r)^4-420 a(r)^3 \sqrt{r} w+90 a(r)^2 r w (-1+8 w)-
$$
$$
-20 a(r) r^{3/2} w (1-12 w+32 w^2)+r^2 w (-2+43 w-192 w^2+256 w^3))]^{-1}
$$
and
$$
A_0=1
$$
The above values of $A_j$ were evaluated with Mathematica on the basis  of the first and third term of (13). Setting all the Taylor expansion coefficients of 
$$
A_{p}z^pF^{(p)}_{p}(z)+A_{p-1}z^{p-1}F^{(p-1)}(z)+\ldots+A_{1}zF^{(1)}(z)+A_0F^{(0)}(z) ,
$$ 
with respect to $K(w)$ to be 0, we arrived at the desired values (we used the relations (4) and (7)). The quantity $C_j$ was evaluated from the first equation of relation (13). The general formula produced by this method for $1/\pi^4$ is 
\[
\]
\begin{equation}
\sum^{\infty}_{n=0}c_{4}(n)(k_rk'_r)^{2n}[A_4n^4+(A_3-6A_4)n^{3}+(A_2-3A_3+11A_4)n^{2}+
$$
$$
+(A_1-A_2+2A_3-6A_4)n^1+A_0]=\frac{g}{\pi^{4}}
\end{equation}
\[
\]
$$
g=105[\pi ^4 (105 a(r)^4-420 a(r)^3 \sqrt{r} w+90 a(r)^2 r w (-1+8 w)-
$$
$$
-20 a(r) r^{3/2} w (1-12 w+32 w^2)+r^2 w (-2+43 w-192 w^2+256 w^3))]^{-1}
$$
with the above values for $A_4$, $A_3$, $A_2$, $A_1$, $A_0$  (page 3),  
where $1-2w=k_r^2$.\\
It is worth pointing out that Ramanujan-type formulas of order $\nu \geq 4$ are presented here for the first time. The  formulas  known previously are of order 1,2,3 (see the references of [28] and Guillera's Pages on the Web).\\ 
\\
\textbf{Application}\\
From Wolfram's Mathworld pages (see Elliptic Lambda Function) and [7] for $r=2$ we have $k_{2}=\sqrt{2}-1$ and $a(2)=\sqrt{2}-1$. Hence we get the formula
\small
$$
\sum^{\infty}_{n=0}c_4(n)\left(-56+40 \sqrt{2}\right)^n[462719+5 \left(292072+56267 \sqrt{2}\right) n+6 \left(268641+81580 \sqrt{2}\right) n^2+
$$
$$
+4 \left(134444+32155 \sqrt{2}\right) n^3-4 \left(36209+34800 \sqrt{2}\right) n^4]
=-\frac{48585495}{\left(-229441+162240 \sqrt{2}\right) \pi ^4}
$$
\normalsize 
where
$$
c_4(n)=\sum^{n}_{s=0}c_2(s)c_2(n-s)\textrm{ and } c_2(n)=\frac{1}{2^{6n}}\sum^{n}_{s=0}\left(^{2s}_{s}\right)^3\left(^{2n-2s}_{n-s}\right)^3 
$$
\\
\textbf{Remark.}\\ 
One can proceed to $k_r$  with higher values of $r\in \bf N\rm$ to obtain more rapidly convergent$1/\pi^4$ series.
\\
\centerline{\textbf{\large A series representing $1/\pi^6$}}
\\
For $\nu=3$ we get\\
\\
\small
\begin{equation}
\sum^{\infty}_{n=0}c_{6}(n)(k_rk'_r)^{2n}[A_6n^6+(A_5-15A_6)n^5+(A_4-10A_5+85A_6)n^4+
$$
$$
+(A_3-6A_4+35A_5-225A_6)n^{3}+(A_2-3A_3+11A_4-50A_5+274A_6)n^{2}+
$$
$$
+(A_1-A_2+2A_3-6A_4+24A_5-120A_6)n^1+A_0]=\frac{g}{\pi^6}
\end{equation}
\\
\normalsize
The coefficients $A_j$ are obtained as in the case $\nu=2$. 
We present some examples\\
\\
\textbf{i)} For $r=2$ is
\[
\]
\small
$$
\sum^{\infty}_{n=1}c_6(n)(-56+40\sqrt{2})^n\times$$
$$\times[1+\frac{\left(28335508172-240070543 \sqrt{2}\right) n}{12623771801}-\frac{\left(-22911684702+3047538900 \sqrt{2}\right) n^2}{12623771801}-
$$
$$
-\frac{\left(-6110502200+5456734120 \sqrt{2}\right) n^3}{12623771801}-\frac{\left(1196112280+3649618320 \sqrt{2}\right) n^4}{12623771801}-
$$
$$
-\frac{\left(505494672+788011092 \sqrt{2}\right) n^5}{12623771801}+\frac{\left(463408744+244639040 \sqrt{2}\right) n^6}{37871315403}]=
$$
$$
=\frac{3465}{\left(629823301-445352320 \sqrt{2}\right) \pi ^6} .
$$
\normalsize
\[
\]
\textbf{ii)} For $r=7$ we have $k^2_7=\frac{8-3\sqrt{7}}{16}$ and $a(7)=\frac{\sqrt{7}-2}{2}$, then\\
\\
\\
\small
$$
\sum^{\infty}_{n=0}\frac{c_6(n)}{64^n}[1+\frac{913150 n}{307323}-\frac{75313 n^2}{102441}-\frac{4998980 n^3}{307323}-\frac{1126755 n^4}{34147}-\frac{1080450 n^5}{34147}-\frac{453789 n^6}{34147}]
$$
$$
=-\frac{14417920}{34147 \pi ^6} .
$$
\\
\normalsize
\textbf{iii)} For $r=15$, $k^2_{15}=\frac{(2-\sqrt{3})^2(\sqrt{5}-\sqrt{3})^2(3-\sqrt{5})^2}{128}$ and $a(15)=\frac{\sqrt{15}-\sqrt{5}-1}{2}$\\
then
\scriptsize
$$
\sum^{\infty}_{n=0}c_6(n)\left(\frac{47-21 \sqrt{5}}{128} \right)^n[1+\frac{\left(2877117109830+924178552332 \sqrt{5}\right) n}{293049243769}+
$$
$$
+\frac{\left(15689590644975+6660423786240 \sqrt{5}\right) n^2}{293049243769}+\frac{\left(51863088153600+23066524139820 \sqrt{5}\right) n^3}{293049243769}+
$$
$$
+\frac{\left(106483989569175+47630637457200 \sqrt{5}\right) n^4}{293049243769}+
$$
$$
+\frac{\left(130261549416750+58266415341540 \sqrt{5}\right) n^5}{293049243769}+\frac{\left(75619648012725+33817435224300 \sqrt{5}\right) n^6}{293049243769}]=
$$
$$
=\frac{20185088}{\left(11556387-5162500 \sqrt{5}\right) \pi ^6}
$$
\[
\]
\normalsize

\section{A reduction formula for the  Elliptic Alpha function $a(25r)$.}

It is clear from the results in  section 2 that to obtain rapidly convergent  series for $1/\pi$ and its even powers one requires values of the alpha function $a(r)$ for large $r\in N$, say $r=5000$ (see also [29] and [11],[25]). In this section we address this problem and conclude with several examples.

From relations (4),(7) and [4] pages 121-122, chapter 21, if we set $y=\pi\sqrt{r}$, $q=e^{-\pi\sqrt{r}}$, $K(k_r)=K[r]$, $k'_r=\sqrt{1-k^2_r}$, then
\begin{equation}
1-24\sum^{\infty}_{n=1}\frac{nq^{2n}}{1-q^{2n}}=\frac{3}{\pi\sqrt{r}}+\left(1+k^2_r-\frac{3\alpha(r)}{\sqrt{r}}\right)\frac{4}{\pi^2}K^2[r]
\end{equation}
From
\begin{equation}
k_{r/4}=\frac{2\sqrt{k_r}}{1+k_r}
\end{equation}
and 
\begin{equation}
\alpha(4r)=(1+k_{4r})^2\alpha(r)-2\sqrt{r}k_{4r} 
\end{equation}
relation (16) becomes
\[
\]
\begin{equation}
1-24\sum^{\infty}_{n=1}\frac{nq^n}{1-q^n}
=\frac{6}{\pi\sqrt{r}}+4\frac{K^2[r]\left(-6 \alpha(r)+\sqrt{r}\left(1+k^2_r\right)\right)}{\pi^2 \sqrt{r}}
\end{equation}
If we set
\begin{equation}
T_{p,r}:=\left(1-24\sum^{\infty}_{n=1}\frac{nq^{2n}}{1-q^{2n}}\right)-p\left(1-24\sum^{\infty}_{n=1}\frac{nq^{2pn}}{1-q^{2pn}}\right)
\end{equation}
then\\
\\
\textbf{Proposition 1.}
\begin{equation}
\frac{1}{m^2_p\sqrt{r}}\alpha(p^2r)=-\frac{(1+k^2_r)}{3}+\frac{p(1+k^2_{p^2r})}{3m^{2}_{p}}+\frac{\pi^2 T_{p,r}}{12K^2[r]}+\frac{\alpha(r)}{\sqrt{r}} .
\end{equation}
\\
The above Proposition 1 relates results of Ramanujan in [4] chapter 21 with the evaluation of  the alpha function and the  evaluations of $\pi$. Solving (21) for $T_{p,r}$ we have
\begin{equation}
T_{p,r}=\frac{4K[r]^2}{\pi^2\sqrt{r}m^2_p}\left[3\alpha(p^2r)-p\sqrt{r}(1+k^2_{p^2r})+\left(-3\alpha(r)+\sqrt{r}(1+k^2_r)\right)m^2_p\right]
\end{equation}
Hence from (19),(20) and (22) we get another interesting formula 
\begin{equation}
1-24\sum^{\infty}_{n=1}\frac{nq^{2pn}}{1-q^{2pn}}=
\frac{3}{\pi\sqrt{r}p}+\frac{4K[r]^2}{\pi^2\sqrt{r}pm^2_p }\left[-3\alpha(p^2r)+p\sqrt{r}(1+k^2_{p^2r})\right]
\end{equation}
where
\begin{equation}
m_p=\frac{K[r]}{K[p^2r]}
\end{equation}
From [4] page 463 Entry 4 we have
\begin{equation}
T_{5,r}=\frac{(x^2+22\cdot q^2 x y +125\cdot q^4y^2)^{1/2}}{(xy)^{1/6}} , 
\end{equation}
where $x=f^6(-q^2)$ and $y=f^6(-q^{10})$.\\
From  Ramanujan's identity
\begin{equation}
A_r:=\frac{1}{R^5(q^2)}-11-R^5(q^2)=\frac{x}{q^2 y}
\end{equation}
where $R(q)$ is the Rogers-Ramanujan continued fraction 
$$R(q)=\frac{q^{1/5}}{1+}\frac{q}{1+}\frac{q^2}{1+}\frac{q^3}{1+}\ldots $$
and 
\begin{equation}
f(-q^2)^6=\frac{2k_rk'_rK[r]^3}{\pi^3q^{1/2}}
\end{equation}
we get
\begin{equation}
T_{5,r}=4\frac{(k_r k'_r)^{2/3}\sqrt{125+22 A_r+A^2_r}}{6\cdot 2^{1/3}A^{5/6}_r}=2^{2/3}\frac{(k_r k'_r)^{2/3}\left[R^{-5}(q^2)+R^5(q^2)\right]}{3\left[R^{-5}(q^2)-11-R^5(q^2)\right]^{5/6}}
\end{equation}
and hence the evaluation
\begin{equation}
-4-24\sum^{\infty}_{n=1}\frac{nq^{2n}}{1-q^{2n}}+120\sum^{\infty}_{n=1}\frac{nq^{10n}}{1-q^{10n}}
=2^{2/3}\frac{(k_r k'_r)^{2/3}\left[R^{-5}(q^2)+R^5(q^2)\right]}{3\left[R^{-5}(q^2)-11-R^5(q^2)\right]^{5/6}}
\end{equation}
But for the evaluation of the Rogers-Ramanujan continued fraction from [26] we have:\\
\\
\textbf{Proposition 2.}(see [30])\\
If $q=e^{-\pi\sqrt{r}}$ and $r$ real positive, then
\begin{equation}
A_r=a_{4r}=\left(\frac{k_rk'_r}{w_rw'_r}\right)^2\left(\frac{w_r}{k_r}+\frac{w'_r}{k'_r}-\frac{w_rw'_r}{k_rk'_r}\right)^3
\end{equation}
with
\begin{equation}
m_5=\frac{w_r}{k_r}+\frac{w'_r}{k'_r}-\frac{w_rw'_r}{k_rk'_r} \textrm{ and } w_r=\sqrt{k_rk_{25r}} \textrm{ , } w'_r=\sqrt{k'_rk'_{25r}}
\end{equation}
Hence we get\\
\\
\textbf{Proposition 3.}
$$
\frac{3\alpha(25r)}{m^2_5\sqrt{r}}-\frac{3\alpha(r)}{\sqrt{r}}=
$$
\begin{equation}
=\frac{5}{m_5^2}(1+k^2_{25r})-(1+k^2_r)-2^{2/3}A_r^{-5/6}(k_rk'_r)^{2/3}[R^{5}(q^2)+R^{-5}(q^2)]
\end{equation}
\textbf{Proof.}\\
From (21), (28) is
$$
\frac{3\alpha(25r)}{m^2_5\sqrt{r}}-\frac{3\alpha(r)}{\sqrt{r}}=
$$
\begin{equation}
=-1+\frac{5}{m_5^2}+\frac{5 k_{25r}^2}{m_5^2}-k_r^2-\frac{2^{2/3}\sqrt{125+22 A_r+A_r^2} (k_rk'_r)^{2/3}}{A_r^{5/6}}
\end{equation}
with 
\begin{equation}
A_{r/4}=\left(\frac{k'_r}{k'_{25r}}\right)^2\sqrt{\frac{k_r}{k_{25r}}}m_5^{-3}
\end{equation}
The multiplier $m_5$ is that of (31) and satisfies also the equation
\begin{equation}
(5m_5-1)^5(1-m_5)=256k^2_rk'^2_rm_5 .
\end{equation}
\\
The next proposition is a conjecture which is most compactly expressed in terms of the quantity  $$Y_{\sqrt{-r}}=\frac{1}{6}\left(R\left(e^{-2\pi  \sqrt{r}}\right)^{-5}-11-R\left(e^{-2\pi \sqrt{r}}\right)^5\right)=\frac{A_r}{6}\eqno{:(a)}
$$
The function $j_r$ which appears is the $j$-invariant (see [6],[17]). For more properties of $j_r$ and $A_r$  see [26]:\\
\\
\textbf{Proposition 4.}\\
As indicated by numerical results,   whenever $(r,5)=$GCD$(r,5)=1$,\\then $deg\left(Y_{\sqrt{-r/5}}\right)=deg\left(j_{\sqrt{-r/5}}\right)$.\\ 
For a given $r\in \bf N\rm$ and $deg\left(Y_{\sqrt{-r/5}}\right)=2$, $4$, or $8$, if the smallest nested root of $j_{\sqrt{-r/5}}$ is $\sqrt{d}$ then we can evaluate the Rogers-Ramanujan continued fraction with integer parameters.\\ 
\textbf{i)} In the case $deg\left(Y_{\sqrt{-r/5}}\right)=2$ then
$$Y_{\sqrt{-r/5}}=\frac{l+m\sqrt{d}}{t}$$
with $$l^2-m^2d=1$$
\textbf{ii)} In the case $deg\left(Y_{\sqrt{-r/5}}\right)=4$ we have\\
\textbf{a)} If $U\neq\frac{125}{64}$, then
\begin{equation}
Y_{\sqrt{-r/5}}=\frac{5}{8}\sqrt{\left(a_0+\sqrt{-1+a_0^2}\right)}\left(\sqrt{5+p}-\sqrt{p}\right) 
\end{equation}
where 
\begin{equation}
Y_{\sqrt{-r/5}}Y^{*}_{\sqrt{-r/5}}=\frac{125}{64}\left(a_0+\sqrt{a_0^2-1}\right) , 
\end{equation}
where $a_0$ is the positive integer-solution of $l^2-m^2d=1$. Hence $l=a_0$ and $m=d^{-1/2}\sqrt{a_0^2-1}$ is a positive integer. The parameter $p$ is positive rational and can be found directly from the numerical value of $Y_{\sqrt{-r/5}}$.\\
\textbf{b)} If $U=\frac{125}{64}$, then
\begin{equation}
Y_{\sqrt{-r/5}}=A+\frac{1}{8}\sqrt{-125+64A^2} ,
\end{equation}
where we set $A=k+l\sqrt{d}$. Hence a starting point for the evaluation of the integers $k$, $l$ is the relation $$l^2=\frac{(A-k)^2}{d}=\textrm{ square of an integer }$$
\textbf{iii)} If $deg\left(Y_{\sqrt{-r4^{-1}5^{-1}}}\right)=4$, then we can evaluate $Y_{\sqrt{-r5^{-1}}}$.\\
If $deg\left(Y_{\sqrt{-r5^{-1}}}\right)=8$, the minimal polynomial of $Y_{\sqrt{-r5^{-1}}}/Y_{\sqrt{-r4^{-1}5^{-1}}}$ is of degree 4 or 8 and is symmetric. Hence it can be reduced to a polynomial of degree at most 4, and hence it is solvable. Thus it remains to evaluatef $Y_{\sqrt{-r4^{-1}5^{-1}}}$, which can be done with the help of step (ii). 
\begin{equation}
Y_{\sqrt{-r5^{-1}}}=\frac{5}{8}\sqrt{a_0+\sqrt{-1+a^2_0}}\left(\sqrt{p+5}-\sqrt{p}\right)2^{-1}\left(\sqrt{x+4}-\sqrt{x}\right) 
\end{equation}
where $x=a_1+b_1\sqrt{d}+c\sqrt{a_2+b_2\sqrt{d}}$, $a_1$, $b_1$, $a_2$, $b_2$, $c$ integers and
$$ Y_{\sqrt{-r5^{-1}4^{-1}}}=\frac{5}{8}\sqrt{a_0+\sqrt{-1+a^2_0}}\left(\sqrt{p+5}-\sqrt{p}\right)
$$
\[
\] 
We give some values of $Y_{\sqrt{-r/5}}=8^{-1}A_{r/5}$.
\begin{equation}
Y_{\sqrt{-1/5}}=\frac{5\sqrt{5}}{8}  
\end{equation}
\begin{equation}
Y_{\sqrt{-2/5}}=\frac{5}{8}\left(5+2\sqrt{5}\right)  
\end{equation}
\begin{equation}
Y_{\sqrt{-3/5}}=\frac{5}{16} \left(25+11 \sqrt{5}\right)  
\end{equation}
\begin{equation}
Y_{\sqrt{-4/5}}=\frac{5}{16} \left(25+13 \sqrt{5}+5 \sqrt{58+26 \sqrt{5}}\right)  
\end{equation}
\begin{equation}
Y_{\sqrt{-5/5}}=\frac{125}{8} \left(2+\sqrt{5}\right)  
\end{equation}
\begin{equation}
Y_{\sqrt{-6/5}}=\frac{5}{8} \left(50+35 \sqrt{2}+3 \sqrt{5 \left(99+70 \sqrt{2}\right)}\right)  
\end{equation}
\begin{equation}
Y_{\sqrt{-9/5}}=\frac{5}{8} \left(225+104 \sqrt{5}+10 \sqrt{1047+468 \sqrt{5}}\right)  
\end{equation}
\begin{equation}
Y_{\sqrt{-12/5}}=\frac{5}{16} \left(1690+975 \sqrt{3}+29 \sqrt{6755+3900 \sqrt{3}}\right)  
\end{equation}
\begin{equation}
Y_{\sqrt{-14/5}}=\frac{5}{8} \left(1850+585 \sqrt{10}+7 \sqrt{5 \left(27379+8658 \sqrt{10}\right)}\right)
\end{equation}
\begin{equation}
Y_{\sqrt{-17/5}}=\frac{5}{8}\left(5360+585 \sqrt{85}+4 \sqrt{3613670+391950 \sqrt{85}}\right)
\end{equation}
\[
\]
\textbf{Example.}\\
For $r=68=4\cdot17$ and from (49) we have $d=85$
$$
x=a_1+b_1\sqrt{85}+c\sqrt{a_2+b_2\sqrt{85}}
$$ 
$$
Y_{\sqrt{-68/5}}/Y_{\sqrt{-17/5}}=2^{-1}\left(\sqrt{x+4}-\sqrt{x}\right)
$$
$$
a_1=2891581250 \textrm{, } b_1=313636050 \textrm{, } c=12960
$$
$$ 
a_2=99557521554 \textrm{, } b_2=10798529365  
$$
hence 
$$
Y_{\sqrt{-68/5}}=Y_{\sqrt{-17/5}}2^{-1}\left(\sqrt{x+4}-\sqrt{x}\right)=
$$
$$
=\frac{5}{16}\left(5360+585 \sqrt{85}+4\sqrt{3613670+391950 \sqrt{85}}\right)\left(\sqrt{x+4}-\sqrt{x}\right)
$$

\[
\]

\centerline{\bf References}\vskip .2in

[1]: M.Abramowitz and I.A.Stegun: Handbook of Mathematical Functions. Dover Publications. (1972)

[2]: B.C.Berndt: Ramanujan`s Notebooks Part I. Springer Verlag, New York. (1985)

[3]: B.C.Berndt: Ramanujan`s Notebooks Part II. Springer Verlag, New York. (1989)

[4]: B.C.Berndt: Ramanujan`s Notebooks Part III. Springer Verlag, New York. (1991)

[5]: B.C. Berndt: Ramanujan's Notebooks Part V. Springer Verlag, New York, Inc. (1998) 

[6]: Bruce C. Berndt and Heng Huat Chan: Ramanujan and the Modular j-Invariant. Canad. Math. Bull. Vol.42(4), (1999). pp.427-440.

[7]: J.M. Borwein and P.B. Borwein: Pi and the AGM. John Wiley and Sons, Inc. New York, Chichester, Brisbane, Toronto, Singapore. (1987) 

[8]: I.S. Gradshteyn and I.M. Ryzhik: Table of Integrals, Series and Products. Academic Press. (1980)
 
[9]: E.T. Whittaker and G.N. Watson: A course on Modern Analysis. Cambridge U.P. (1927)

[10]: I.J. Zucker: The summation of series of hyperbolic functions. SIAM J. Math. Ana.10.192. (1979)

[11]: Bruce.C. Berndt and Heng Huat Chan: Eisenstein Series and Approximations to $\pi$. Page stored in the Web.

[12]: S. Ramanujan: Modular equations and approximations to $\pi$. Quart. J. Math.(Oxford). 45, 350-372. (1914).

[13]: S. Chowla: Series for $1/K$ and $1/K^2$. J. Lond. Math. Soc. 3, 9-12. (1928)

[14]: N.D. Baruah, B.C. Berndt and H.H. Chan: Ramanujan's series for $1/\pi$: A survey. American Mathematical Monthly 116, 567-587. (2009)

[15]: T. Apostol: Modular Functions and Dirichlet Series in Number Theory. Springer

[16]: Bruce.C. Berndt, S. Bhargava and F.G. Garvan: Ramanujan's Theories of Elliptic Functions to Alternative Bases. Transactions of the American Mathematical Society. 347, 4163-4244. (1995)

[17]: D. Broadhurst: Solutions by radicals at Singular Values $k_N$ from New Class Invariants for $N\equiv3\;\; mod\;\; 8$'. arXiv:0807.2976 (math-phy)

[18]: J.V. Armitage W.F. Eberlein: Elliptic Functions. Cambridge University Press. (2006)

[19]: N.D. Baruah, B.C. Berndt: Eisenstein series and Ramanujan-type series for $1/\pi$. Ramanujan J.23. (2010) 17-44

[20] N.D. Baruah, B.C. Berndt: Ramanujan series for $1/\pi$ arising from his cubic and quartic theories of elliptic functions. J. Math. Anal. Appl. 341. (2008) 357-371 

[21]: B.C. Berndt: Ramanujan's theory of Theta-functions. In Theta functions: from the classical to the modern Editor: Maruti Ram Murty, American Mathematical Society. 1993 

[22]: J.M. Borwein and P.B. Borwein: A cubic counterpart of Jacobi's identity and the AGM. Transactions of the American Mathematical Society, 323, No.2, (Feb 1991), 691-701 

[23]: Habib Muzaffar and Kenneth S. Williams: Evaluation of Complete Elliptic Integrals of the first kind at Singular Moduli. Taiwanese Journal of Mathematics, Vol. 10, No. 6, pp 1633-1660, December 2006 

[24]: Bruce C. Berndt and Aa Ja Yee: Ramanujans Contributions to Eisenstein Series, Especially in his Lost Notebook. (page stored in the Web). 

[25]: Nikos Bagis: Eisenstein Series, Alternative Modular Bases and Approximations of $1/\pi$. arXiv:1011.3496. (2010)

[26]: Nikos Bagis: On a General Polynomial Equation Solved by Elliptic Functions. arXiv:1111.6023v1. (2011)  

[27]: Nikos Bagis: Evaluation of Fifth Degree Elliptic Singular Moduli. arXiv:1202.6246v1. (2012)

[28]: Wadim Zudilin: Ramanujan-type formulae for $1/\pi$: A second wind?$^*$. arXiv:0712.1332v2. (2008)

[29]: Nikos Bagis: Ramanujan type $1/\pi$ Approximation Formulas.\\arXiv:1111.3139v1. (2011) 

[30]: Nikos Bagis. 'Parametric Evaluations of the Rogers Ramanujan Continued Fraction'. International Journal of Mathematics and Mathematical Sciences. Vol. 2011

\end{document}